\newtheorem{lemma}{Lemma}[section]
\newtheorem{prop}[lemma]{Proposition}
\newtheorem{theorem}[lemma]{Theorem}
\newtheorem{cor}[lemma]{Corollary}
\newtheorem{rem}[lemma]{Remark}

\newcommand{\kla}{\left ( }
\newcommand{\mer}{\right ) }

\newcommand{\nz}{{\rm  I\! N}}
\newcommand{\nen}{n \in \nz}
\newcommand{\rz}{{\rm  I\! R}}
\newcommand{\kz}{{\rm  I\! K}}

\newcommand{\dop}{{\rm I\! D}{\cal L}}

\newcommand{\p}{\hspace{.05cm}}
\newcommand{\pl}{\hspace{.1cm}}
\newcommand{\pll}{\hspace{.3cm}}
\newcommand{\pla}{\hspace{1.5cm}}
\newcommand{\hz}{\vspace{0.5cm}}

\newcommand{\Om}{\Omega}
\newcommand{\om}{\omega}
\newcommand{\al}{\alpha}
\newcommand{\be}{\beta}
\newcommand{\si}{\sigma}

\newcommand{\vare}{\varepsilon}
\newcommand{\ds}{D_{\si}}

\newcommand{\lzn}{\ell_2^n}

\newcommand{\lif}{\ell_{\infty}}
\newcommand{\lih}{\ell_{\infty,\infty,1/2}}
\newcommand{\ce}{{\tt c}_o}
\newcommand{\com}{{\cal L}(}

\newcommand{\pxe}{\pi_{X,1}}
\newcommand{\pxz}{\pi_{X,2}}

\newcommand{\noo}{\left \|}
\newcommand{\rrm}{\right \|}

\newcommand{\bet}{\left |}
\newcommand{\rag}{\right |}
\newcommand{\summ}{\sum\limits}

\documentstyle{article}
\begin{document}

\title{Comparing gaussian and Rademacher cotype for operators
 on the space of continous functions}

\author{Marius Junge}

\date{}

\maketitle

\begin{abstract} We will prove an abstract comparision principle which
translates gaussian cotype in Rademacher cotype conditions and vice
versa. More
precisely, let $2\!<\!q\!<\!\infty$ and $T:\,C(K)\,\to\,F$ a linear,
continous
operator.

\begin{enumerate}
\item T is of gaussian cotype q if and only if

    \[ \kla \summ_1^n \kla \frac{\noo Tx_k
    \rrm_F}{\sqrt{\log(k+1)}}\mer^q
     \mer^{1/q} \, \le \pll c \pll
    \noo \summ_1^n \vare_k x_k \rrm_{L_2(C(K))} , \]

for all sequences with $(\noo Tx_k \rrm)_1^n$ decreasing.
\item T is of Rademacher cotype q if and only if

     \[ \kla \summ_1^n \kla \noo Tx_k \rrm_F \,\sqrt{\log(k+1)}\,\mer^q
     \mer^{1/q} \, \le \pll c \pll
    \noo \summ_1^n g_k x_k \rrm_{L_2(C(K))} , \]

for all sequences with $(\noo Tx_k \rrm)_1^n$ decreasing.
\end{enumerate}

Our methods allows a restriction to a fixed number of
vectors and complements the corresponding results of Talagrand.
\end{abstract}

\setcounter{section}{0}


\section*{Introduction}

One problem in the local theory of Banach spaces consits in the
description
of Rademacher cotype and gaussian cotype for operators on
$C(K)$-spaces.
A quite satisfactory answer for the Rademacher cotype was given by
Maurey.
He connected cotype conditions with summing conditions (see
\cite{MAU}):

\hz
\begin{theorem}{\rm [Maurey]} \label{MAU}  Let $2\!<\!q\!<\!\infty$ and
$T:C(K)\,\to\,F$.
 Then the following are eqiuvalent:
 \begin{enumerate}
  \item T is absolutely $(q,2)$-summing, i.e. for all
  $(x_k)_{k\in\nz}\,
	\subset C(K)$ one has
    \[ \kla \summ_k \noo Tx_k \rrm^q \mer^{1/q} \, \le \, c_0 \,
      \sup_{t\in K} \kla \summ_k \bet x_k(t) \rag^2 \mer^{1/2}\pll . \]

  \item T has Rademacher cotype q, i.e. for all $(x_k)_{k\in\nz}\,
	\subset C(K)$ one has
    \[ \kla \summ_k \noo Tx_k \rrm^q \mer^{1/q} \, \le \, c_0 \,
     \noo \summ_k \vare_k x_k \rrm_{L_2(C(K))}       \pll . \]

  \item T is absolutely $(q,1)$-summing, i.e. for all
  $(x_k)_{k\in\nz}\,
	\subset C(K)$ one has
    \[ \kla \summ_k \noo Tx_k \rrm^q \mer^{1/q} \, \le \, c_0 \,
    \sup_{t\in K} \summ_k \bet x_k(t) \rag \pll . \]
 \end{enumerate}
\end{theorem}
\hz

Later on, Pisier gave another approach to this type of results via
factorization theorems. This way was pursued by Montgomery-Smith,
\cite{MSM}, and Talagrand, \cite{TAL}, to give a characterization
of gaussian cotype q.

\hz
\begin{theorem}{\rm [Talagrand]} \label{TA} Let $2\!<\!q\!<\!\infty$
and
 $T:C(K)\,\to \,F$. Then the following are equivalent.

\begin{enumerate}
   \item T has gaussian cotype q, i.e. for all $(x_k)_{k\in\nz}\,
	\subset C(K)$ one has
    \[ \kla \summ_k \noo Tx_k \rrm^q \mer^{1/q} \, \le \, c_1 \,
     \noo \summ_k g_k x_k \rrm_{L_2(C(K))}       \pll . \]

   \item T satisfies the following summing condition, i.e.$\p$for all
   $(x_k)_{k\in\nz}\,
\subset C(K)$ such that $(\noo Tx_k\rrm)_1^n$ is decreasing one has
    \[ \kla \summ_k \kla \frac{\noo Tx_k \rrm}{\sqrt{\log(k+1)}}\mer^q
     \mer^{1/q} \, \le \, c _2\,
    \sup_{t\in K} \summ_k \bet x_k(t) \rag \pll . \]

  \item T factors through an Orlicz space $L_{t^q(\log
  t)^{q/2},\,1}(\mu)$ for
     some probability measure $\mu$ on K.
\end{enumerate}
\end{theorem} \hz

The main new ingredient of this theorem is a factorization theorem for
gaussian processes derived from the existence of majorizing measures,
see \cite{TA1}.

\hz

We will give a more abstract approach to gaussian cotype conditions
which can be considered as a complement to Talagrand's results.
Independently of him we discovered the connection between gaussian
cotype
and summing properties with the modified
$\ell_q$ space in condition 2 of theorem 2.  In order to be precise,
let us give the following definition. For a maximal, symmetric sequence
space
X and $T:\,E\,\to\,F$ we define

\begin{eqnarray*}
\pi_{X,q}^n(T) &:=& \sup\left\{\,
\noo \summ_1^n \noo Tx_k \rrm_{F} \, e_k \rrm_X \pll \bet \pll
  \sup_{a \in B_{E^*}}
  \, \kla \summ_1^n \bet <x_k,a> \rag^q
  \mer^{1/q}\,\le\,1\,\right\}\right.
    \,,\\
rc_X^n(T) &:=&\sup\left\{\, \noo \summ_1^n \noo Tx_k\rrm_F \, e_k
\rrm_X \pll
\bet\pll
  \noo \summ_1^n \vare_k x_k
  \rrm_{L_2(E)}\,\le\,1\,\right\}\right.\,,\\
gc_X^n(T) &:=& \sup\left \{\, \noo \summ_1^n \noo Tx_k\rrm_F \, e_k
\rrm_X \pll
\bet\pll
  \noo \summ_1^n g_k x_k \rrm_{L_2(E)} \,\le\,1\,\right\} \right.\,.\\
\end{eqnarray*}

An operator is said to be (absolutely) $(X,q)-summing$, of $Rademacher$
$cotype$ $X$, of $gaussian$ $cotype$ $X$ if
$\pi_{X,q}\,:=\,\sup_{\nen}\,\pi_{X,q}^n$,
$rc_X\,:=\,\sup_{\nen}\,rc_X^n$,
$gc_X\,:=\,\sup_{\nen}\,gc_X^n$
is finite, respectively. In contrast to
Talagrand we follow Maurey's approach and prove \hz

\begin{theorem} \label{genmau} Let $2\!<\!q\!<\!\infty$, X a q-convex,
  maximal, symmetric sequence space and  $T:C(K)\,\to\,F$. Then the
  following are equivalent:

 \begin{enumerate}
  \item T is $(X,2)$-summing.
  \item T is of Rademacher cotype X.
  \item T is $(X,1)$-summing.
 \end{enumerate}

Furthermore, there exists a constant c only depending on q and X such
that
\[ \pxz^n(T) \: \le \,c\, \pxe^n(T) \pll . \]

\end{theorem}
\hz

The main idea for the proof of the theorem above is a reduction to
Maurey's
result via quotient formulas. These formulas are contained in chapter 2
and
have already be seen to be helpful in the theorey of summing operators.
Their
proof goes back to a joint work of Martin Defant and the author, see
\cite{DJ}.
The comparision principle between gaussian and Rademacher cotype for
operators on $C(K)$-spaces is formulated in \hz

\begin{theorem} \label{RADGAU} Let $2\!<\!q\!<\!\infty$, X a q-convex,
 maximal, symmetric sequence space. If Y denotes the space of diagonal
 operators
between $\ell_{\infty,\infty,1/2}$ and $X$ one has for all operators
$T:C(K) \to F$ and $n\in\nz$

 \[\frac{1}{c} rc_Y^n(T) \, \le gc_X^n(T) \,\le c\, rc_Y^n(T) \pll,\]

 where c is a constant depending on q and X only.

\end{theorem} \hz

The philosophy is quite simple. The difference between gaussian and
Rademacher cotype has to be corrected in the summing property with the
factor $\sqrt{\log(k+1)}$. This becomes clear if we apply this first
for
the space $X\,=\,\ell_q$. Then we see that an opertor
$T:\,C(K)\,\to\,F$
is of gaussian cotype q if and only if

    \[ \kla \summ_k \kla \frac{\noo Tx_k
    \rrm_F}{\sqrt{\log(k+1)}}\mer^q
     \mer^{1/q} \, \le \, c \,
    \noo \summ_1^n \vare_k x_k \rrm_{L_2(C(K)} , \]

for all sequences with $(\noo Tx_k \rrm)_1^n$ decreasing. Applying the
result
for $Y\,=\,\ell_q$ we see that T is of Rademacher cotype q if and only
if

     \[ \kla \summ_k \kla \noo Tx_k \rrm_F \,\sqrt{\log(k+1)}\,\mer^q
     \mer^{1/q} \, \le \, c \,
    \noo \summ_1^n g_k x_k \rrm_{L_2(C(K)} , \]

for all sequences with $(\noo Tx_k \rrm)_1^n$ decreasing. Let us also
note
that our approach enables us to fix the number of vectors in
consideration.
For example, this restriction to n vectors can be used to prove that
for an
opertor of rank n the gaussian cotype q-norm is attained on n disjoint
functions in $C(K)$. Another application is given in the study of weak
cotype
operators. \hz

\section*{Preliminaries}

We use standard Banach space notations. In particular, $c_0$, $c_1$, ..
will
denote different absolute constants and they can vary whithin the
text.
The symbols X, Y, Z are reserved for sequence spaces. Standard
references on
sequence spaces and Banach lattices are the monograph of Lindenstrauss
and
Tzafriri, \cite{LTI,LTII}. The symbols $E$, $F$ will always denote
Banach sapces
with unit balls $B_E$, $B_F$ and duals $E^*$, $F^*$. Basic information
on operator ideals and s-numbers can be found in the monograph of
Pietsch,
\cite{PIE}. The ideal of linear operators is denoted by ${\cal L}$.
\hz

The classical sequence spaces $\ce$, $\ell_p$ and $\ell_p^n$,
$1\!\le\!p\!\le\!\infty$, $\nen$ are defined in the usual way. From
the context it will be clear whether we mean the space $\ce$ or the
absolut constant $c_0$. A generalization of the classical $\ell_p$
spaces
is the class of Lorentz-Marcinkiewicz spaces. For a given continous
function $f:\nz \to \rz_{>0}$ with $f(1)=1$ the following two indices
are defined

\[ \al_f \, := \, \inf \{\, \al \, | \, \exists M<\infty \, \forall
   t,s\ge 1: \: f(ts)\, \le \,M t^{\al} f(s) \, \} \pll, \]

\[ \be_{f} \, := \, \sup \{\, \be \, | \, \exists c>0 \, \forall
   t,s\ge 1: \: f(ts)\, \ge \,c t^{\be} f(s) \, \} \pll. \]

These two indices play an important r\^ole in the study of the space
$\ell_{f,q}$, $1\!\le\!q\!\le\!\infty$ consisting  of all sequences
$\si \in \lif$ such that

\[\noo \si \rrm_{f,q} \pl
  := \pl \kla \summ_n \kla f(n)\, \si_n^* \mer^q n^{-1} \mer^{1/q}
  \,<\,
  \infty \, \pll.   \]

For $q\!=\!\infty$ the needed modification is given by

\[ \noo \si \rrm_{f,\infty}  \pl := \pl  \sup_{\nen} \, f(n)\,\si_n^*
   \pl <\pl \infty . \]

Here and in the following $\si^*\,=\,(\si_n^*)_{\nen}$ denotes the
non-increasing rearrangement of $\si$.
\hz

In the introduction the notions of $(X,q)$-summing, Rademacher cotype X
and
gaussian cotype X are already defined. If $X = \ell_p$ we will shortly
speak of $(p,q)$-summing opertors or norms, Rademacher cotype p, etc.
(possibly restricted to n vectors). In this context it is convenient to
use
an abbreviation for the right hand side of the definition of summing
operators. For a sequence $(x_k)_1^n$ in a Banach space E we write

\[ \om_q(x_k)_1^n \pl := \pl
 \sup_{a \in B_{E^*}} \, \kla \summ_1^n \bet <x_k,a> \rag^q
 \mer^{1/q}\pll.\]

Let us note that this expression coincides with the operator norm of

\[ u \, := \, \summ_1^n e_k \otimes x_k \pl \in \com \ell_{q'}^n,E)
\pll,\]

where $q'$ is the conjugate index of $q$ satisfying $\frac{1}{q}\,+\,
 \frac{1}{q'}\,=\,1$.\hz

In the following $(\vare_n)_{\nen}$, $(g_n){\nen}$ will denote a
sequence of
independent normalized Bernoulli ($Rade- macher$) variables or
$gaussian$
variables respectively. They are defined on a probability space
$(\Om,\mu)$.
Here Bernoulli variable means

\[ \mu(\vare_n \, =\, +1) \pl = \pl \mu(\vare_n \,= \, -1)
\pl =\pl \frac{1}{2} \pl .\]

A very deep result in the theory of gaussian processes is Talagrand's
factorization theorem, see \cite{TA1}.
\hz

{\it
$\bf (*)$ \pll  There is an absolut constant $c_1$ such that for all
sequence $(x_k)_1^n \in C(K)$ with

\[\noo \summ_1^n g_k x_k \rrm_{L_2(X)} \: \le \: 1 \pll . \]

\hspace{0.8 cm} there are operators $u: \lzn \to \ce$, $R:\ce \to C(K)$
 with $\noo u \rrm$ $\noo R \rrm \le c_1$ such that

\[RD_{\sigma}u(e_k) \: =\: x_k \pll,\]

\hspace{0.8 cm} where $D_{\sigma}$ is the diagonal operator with

\[ \sigma_k \:=\: \frac{1}{\sqrt{\,\log(k+1)}}  \pll . \]

}

Finally some s-numbers are needed. For an operator $T\in\com E,F)$ and
$\nen$
the $n$-th $approximation$ $number$ is defined by

\[ a_n(T) \pl :=\pl \inf\{\, \noo T-S \rrm \, | \,rank(S)\,< \, n \,\}
\pla ,\]

whereas the $n$-th $Weyl\,number$ is given by

\[ x_n(T) \pl :=\pl \sup\{\, a_n(Tu)\, |\, u \in \com \ell_2,E) \,
\mbox{with}
 \, \noo u \rrm \,\le \, 1\,\} \pla .\]

\section{Maximal symmetric sequence spaces}

\setcounter{lemma}{0}

In the following we will denote the set of all finite sequences by
$\phi$ and
the sequence of unit vectors in $\ell_{\infty}$ by $(e_k)_k$. For every
sequence $\sigma \!=\! (\sigma_k)_k \subset \ell_{\infty}$, $n \in \nz$
we set
$P_n(\sigma)\,:=\, \summ_1^n \sigma_k e_k$.\hz

A $maximal\, sequence\, space \, (X,\noo \cdot \rrm)$ is a Banach space
satisfying
the following conditions.

\begin{enumerate}
 \item $\ell_1 \subset X \subset \lif$ and $\noo e_k \rrm \, = \, 1$
 for all
     $k \in \nz$.
 \item If $\sigma \in X$ and $\al \in \lif$ then the pointwise product
    $\al \sigma \in X$ with $\noo \al\sigma \rrm \, \le  \, \noo
    \si\rrm_X\,
	 \noo \al \rrm_{\infty}$.
 \item $\sigma \in X$ if and only if $(\,\noo P_n \rrm \, )_n$ is
 bounded and
   in this case
  \[ \noo \si \rrm \,= \, \sup_{n \in \nz} \noo P_n \rrm \pll . \]
\end{enumerate}

For $\nen$ and $\si \, =\,(\si_k)_1^n \, \subset \kz^n$ we set
$\noo \si \rrm := \noo (\si_k)_1^n \rrm := \noo \summ_1^n \si_k e_k
\rrm$. The
sequence dual of X is defined by

\[ X^+ \, := \, \{ \: \tau \in \lif \: | \: \noo \tau \rrm_+ \, := \,
    \sup_{\si \in B_X} \bet \summ_k \si_k \tau_k \rag \, < \, \infty \:
    \}
    \pll . \]

Then $(X^+, \noo \cdot \rrm_+)$ is also a maximal sequence space. We
observe
that $\noo \tau \rrm_{X^*} \, =\, \noo \tau \rrm_+$ holds for all $\tau
\in
\phi$. Thus $X^{++}\,=\,X$ with equal norms. For two maximal sequnce
spaces X, Y
we denote by $\dop(X,Y)$ the space of continous diagonal operators from
X to Y with
the operator norm. A maximal sequence space is $symmetric$ if in
addition $\si \in X$ if and
only  if $\si^* \in X$ with $\noo \si^* \rrm_X \, =\, \noo \si
\rrm_X$.\hz

Essentially for the following is the definition of p-convex sequence
spaces.
Let $1\!\le\!p\!<\!\infty$. A maximal sequence space is $p\!-\!convex$
if
there is a constant $c\!>\!0$ such that for all
$\nen$ and $(x_k)_1^n \subset X$

\[ \noo \kla \summ_1^n \bet x_k \rag^p \mer^{1/p} \rrm \,
  \le c \, \kla \summ_1^n \noo x_k \rrm^p \mer^{1/p} \pll . \]

The best constant c satisfying the above condition will be denoted by
$M^p(X)$.
Obviously, every maximal sequence space is 1-convex. On the other hand
we observe
\[X^+ \, =\, \dop(X,\ell_1) \:\: \mbox{  and thus  }\:\: X
\,=\,\dop(X^+,\ell_1) \,. \]
More generally, one has\hz

\begin{prop}\label{qconvex} Let $1\!\le\!p\!<\!\infty$ and X a maximal
  sequence space. Then the following are equivalent:
 \begin{enumerate}
	\item X is p-convex.
	\item The homogenous expression $\noo\bet\si\rag^{1/p}\rrm_X^p$
	 is equivalent to a norm $\noo\cdot\rrm_p$ with
	 \[\frac{1}{c} \noo \si \rrm_X \, \le \, \noo \bet \si \rag^p
	 \rrm_p^{1/p} \, \le  \, \noo \si \rrm_X \pll. \]
	\item There exists a maximal sequence space Y such that
		\[X \, \cong \, \dop(Y,\ell_p) \pll. \]
  \end{enumerate}
 Moreover, in this case we can choose $Y\,=\,\dop(X,\ell_p)$ and have
   \[\frac{1}{M^p(X)}\, \noo \si \rrm_X \,\le \, \noo D_{\si} \rrm \,
   \le \noo \si\rrm_X\pll .\]
\end{prop}\hz
$\bf Proof:$ The equivalence between 1. and 2. is classical and can be
found
 for example in \cite{LTII}. Now we proof $2. \Rightarrow 3.$ We denote
 by
 $X_p$ the maximal sequence space defined by the norm $\noo \cdot
 \rrm_p$.
 We set $Y:= \dop(X,\ell_p)$. Cleary, we have $X\subset
 \dop(Y,\ell_p)$. By the observations above we have

\begin{eqnarray*}
 \frac{1}{c} \, \noo \si \rrm &\le& \noo \bet\si  \rag^p \rrm_p^{1/p}\\
 &=& \sup_{\tau \in B_{(X_p)^+}} \bet \summ_k \bet\si_k\rag^p  \tau_k
      \rag^{1/p} \\
 &\le& \noo \si \rrm_{\dop(Y,\ell_p)} \sup_{\tau \in B_{(X_p)^+}} \noo
 \bet
  \tau \rag^{1/p} \rrm_{\dop(X,\ell_p)}\\
 &=& \noo \si \rrm_{\dop(Y,\ell_p)} \sup_{\tau \in B_{(X_p)^+}}
     \sup_{\rho \in B_X} \kla \summ_k \bet \tau_k \rag \bet \rho_k
     \rag^p \mer^{1/p} \\
 &=& \noo \si \rrm_{\dop(Y,\ell_p)} \sup_{\rho\in B_X} \noo \bet \rho
      \rag^p \rrm_p^{1/p} \\
 &\le& \noo \si \rrm_{\dop(Y,\ell_p)} \pll .\\
\end{eqnarray*}

For the proof of $3.\Rightarrow 1.$ we can asssume that
$X=\dop(Y,\ell_p)$
with equal norms. The definition of the norm implies for $(x_j)_1^n
\subset X$

\begin{eqnarray*}
  \noo \kla \summ_1^n \bet x_j \rag^p \mer^{1/p} \rrm &=&
    \sup_{\tau \in B_Y} \kla \summ_k \summ_{j=1}^n \bet x_j(k)\rag^p
     \bet \tau_k \rag^p \mer^{1/p} \\
  &=&  \sup_{\tau \in B_Y} \kla \summ_{j=1}^n \summ_k \bet x_j(k)
  \tau_k\rag^p
      \mer^{1/p} \\
  &\le&  \kla \summ_{j=1}^n \sup_{\tau \in B_Y} \kla \summ_k \bet
  x_j(k) \tau_k\rag^p
     \mer \mer^{1/p} \\
  &=& \kla \summ_{j=1}^n \noo x_k \rrm^p \mer^{1/p} \pll . \\[-5ex]
 \end{eqnarray*} \hfill $\Box$

\begin{rem} \label{example} $\bf i)$ An Orlicz sequence space
\[\ell_{\phi} \, := \, \{\, \si \in \lif \, | \, \summ_k \phi(\si_k) \,
< \,
    \infty \,\} \]
is p-convex if and only if $\phi(t\lambda) \, \le\,c
\lambda^p\,\phi(t)$.\hfill
\hz

$\bf ii)$ The criterion above is very useful to study the p-convexity
of a Lorentz-Marcinkiewicz sequence space $\ell_{f,q}$. It was observed
in
\cite{COB} that for $p\!\le\!q$ and $0 \!<\be_{f}\!\le\!\al_f\!< \!1/p$
one has

\[ \noo \bet \si\rag^{1/p} \rrm_{f,q}^p \pll \sim \pll \noo \kla
\frac{1}{n}
    \summ_1^n \si_k^* \mer \rrm_{f^p,q/p} \pll .\]

Since the right hand side is a norm, see again \cite{COB}, the
conditions
above imply the p-convexity of $\ell_{f,q}$.

\end{rem}

\section{Quotient formulas for summing properties}
\setcounter{lemma}{0}

We will start with a quotient formula for $(X,q)$-summing operators.
\hz

\begin{prop}\label{quot1} Let $1\!\le\!r\!\le \!q \!\le\!\infty$,
Y a maximal, symmetric sequence space and $X\cong \dop(Y,\ell_q)$. Then
we
 have for all $\nen$ and $T\in\com E,F)$
\begin{eqnarray*}
 \pi_{X,r}^n(T) &=& \sup\{\,\pi_{q,r}^n(D_{\si}RT)\, |\, R \in \com F,
   \lif),\, \ds \in \com\lif,\lif),\,\mbox{with} \noo R\rrm,
   \noo \si \rrm_Y\le 1 \, \} \, .\\
\end{eqnarray*}
\end{prop}
\hz

$\bf Proof:$ $\bf "\le"$ Let $(x_k)_1^n \subset E$ with.
For $\vare \!>\!0$ there exists a $\si \in B_Y$ with
\[ \noo \,\summ_1^n \noo Tx_k \rrm \, e_k \, \rrm_X \: \le \:
 (1+\vare) \, \kla \summ_1^n \bet \noo Tx_k \rrm \si_k \rag^q
  \mer^{1/q} \: . \]
Let $y_k^* \in B_{F^*}$ with $<y_k^*,Tx_k>\,=\,\noo Tx_k\rrm$. If we
define
$R\,:=\, \summ_1^n y_k^* \otimes e_k \, \in \com F,\lif)$ we obtain
\begin{eqnarray*}
 \frac{1}{1+\vare}\, \noo \,\summ_1^n \noo Tx_k \rrm \, e_k \, \rrm_X
 &\le&
   \kla \summ_1^n \bet <y_k^*,Tx_k> \si_k \rag^q \mer^{1/q} \\
 &\le& \kla \summ_1^n  \sup_j \bet<y_j^*,Tx_k> \si_j \rag^q \mer^{1/q}
 \\
 &\le& \pi_{q,r}^n(\ds R T) \, \om_r(x_k)_1^n \pla .\\[0.3 cm]
\end{eqnarray*}
$\bf"\ge"$ Let $\si \in B_Y$ and $R\in \com F,\lif)$ with $\noo R
  \rrm \, \le \, 1$. By the maximality of $(X,r)$-summing operators
  there is no restriction to assume $R \in \com F, \lif^m)$ for some
 $m \in \nz$. Now we will use a duality argument. Following the proof
 of
 theorem 1. in \cite{DJ} there is an operator $S \in \com \lif^m ,E)$
 with
	\[ \pi_{q,r}^n(\ds RT) \, = \, trace(S\ds RT)\pla \mbox{ and }
	\pla
	  S\,=\,BD_{\tau}P\, , \]
 where $B \in \com \ell_{r'}^n,E)$ with $\noo B\rrm\le1$, $\tau \in
 B_{\ell_{q'}^n}$ and there is an increasing sequence
 $(l_k)_1^n \in \{1,..,m\}$ such that
	\[ P\: = \: \summ_1^n e_{l_k} \otimes e_k \, \in \com
	\lif^m,\lif^n) \pll .\]
 Therefore we deduce
\begin{eqnarray*}
 trace(S\ds RT) &=& trace(D_{\tau} P \ds RTB) \\
 &=& \summ_1^n \tau_k <e_{l_k}\, ,\, \ds RTB(e_k)> \\
 &\le& \summ_1^n \bet \tau_k\,\si_{l_k} \rag  \noo RTB(e_k) \rrm \\
 &\le& \kla \summ_1^n \kla\bet\si_{l_k} \rag  \noo RTB(e_k)\rrm \mer^q
   \mer^{1/q}\\
 &\le& \noo \si \rrm_Y \, \pi_{X,r}^n(RT) \, \noo B \rrm \\[+0.2cm]
 &\le& \pi_{X,r}^n(T) \pla . \\*[-1.2cm]
\end{eqnarray*}\hfill $\Box$\hz

We can now prove the generalized Maurey theorem.\hz

\begin{theorem} \label{gema} Let $1\!\le\!r\!<\!q\!\le\!\infty$, X
 a q-convex maximal, symmetric sequence space and $\nen$. Then for all
 operators $T\in \com C(K),F)$ one has
\[ \pi_{X,r}^n(T) \: \le \: c_0 \, M^q(X) \, \frac{1}{r} \,
    \kla \frac{1}{r}\,-\,\frac{1}{q} \mer^{-1/q'} \, \pxe^n(T)\,. \]
\end{theorem} \hz
$\bf Proof:$ By proposition \ref{qconvex} we can assume that there
exists a
 maximal, symmetric sequence space Y with $X\cong\dop(Y,\ell_q)$. By
 the
 classical Maurey theorem, for the constants see \cite{TJM}, we deduce
 from proposition \ref{quot1}
\begin{eqnarray*}
 \pi_{X,r}^n(T) &\le& \, M^q(X)\,\sup\{\,\pi_{q,r}^n(D_{\si}RT)\, |\,
     R \in \com F,\lif),\, \ds \in \com\lif,\lif),\,\mbox{with} \noo
     R\rrm,
   \noo \si \rrm_Y\le 1 \, \} \\
  &\le& \, M^q(X) \, c_0 \, \frac{1}{r} \,\kla
  \frac{1}{r}\,-\,\frac{1}{q}
    \mer^{-1/q'}\, \times \\
 &  & \pla  \sup\{\,\pi_{q,1}^n(D_{\si}RT)\, |\, R \in \com F,\lif),\,
    \ds \in \com\lif,\lif),\,\mbox{with} \noo R\rrm,
    \noo \si \rrm_Y\le 1 \, \}  \\[0.3cm]
  &=& \, c_0 \, M^q(X) \, \pxe^n(T) \pl .\\[-1.3cm]
\end{eqnarray*} \hfill $\Box$\hz

\begin{rem}\label{intthm} Now it is again well-known, see \cite{MAU},
how
 to derive from the above theorem the equivalence between Rademacher
 cotype
 conditions and summing properties as stated in the introduction as
 theorem 3, namely
  \[ \pxe^n(T) \, \le \, rc_X^n(T) \, \le \, \sqrt{2} \, \pxz^n(T)
\ \le \, c_0 \, M^q(X) \kla \frac{1}{2}\,-\,\frac{1}{q} \mer^{-1/q'}
    \, \pxe^n(T)\,. \]
\end{rem} \hz

At the end of this chapter we will prove another quotient formula which
is
more adapted for operators on $C(K)$-spaces. \hz

\begin{prop} \label{quot2} Let $Y$,$Z$ be maximal, symmetric sequence
spaces
 and $X\, =\, \dop(Y,Z)$. then we have for all $T \in \com E,F)$ and
 $\nen$

 \begin{eqnarray*}
   \pi_{X,1}^n(T) &=& \sup\{\,\pi_{Z,1}^n(TRD_{\si})\, |\, R \in \com
   \lif,
    E),\, \ds \in \com\lif,\lif),\,\mbox{with} \noo R\rrm,
    \noo \si \rrm_Y\le 1 \, \} \, .\\
 \end{eqnarray*}
\end{prop}
\hz
$\bf Proof:$ $\bf "\le"$ can be proved exactly as in proposition
\ref{quot1}.
\hfill

$\bf "\ge"$ Again by maximality we can assume $R \in \com \lif^m ,E)$
and $\ds
  \in \com\lif^m,\lif^m)$ with $\noo R\rrm, \noo \si \rrm_Y \, \le \,
  1$. We
  have to show that for all $S\in \com\lif^n,\lif^m)$ with $\noo S \rrm
  \,
  \le\, 1$ we have

\[ \noo \,\summ_1^n \noo TR\ds S(e_k) \rrm_F \, e_k \, \rrm_Z \: \le
\:
     \pxe^n(T) \,. \]

  By a lemma of Maurey, calculating essentially the extreme points
  of operators from $\lif^n$ to $\lif^m$, see \cite{MAU}, and using
  the convexity of $Z$ we can assume that S has the form

\[S\,=\, \summ_1^n e_k \otimes g^k \pla . \]

  Here the $(g^k)$'s have disjoint support and satisfy $0\!<\!\noo
   g^k \rrm_{\lif^m}\! \le\! 1$. Now we define

\[ J \, :=\,R \,\kla \summ_1^n e_k \otimes \frac{\ds g^k}{\noo \ds g^k
   \rrm_{\infty}} \mer \: \in \com \lif^n,E) \]

  and  $\tau \,:=\,\kla\noo \ds g^k \rrm_{\infty}\mer_1^n$. We observe
  that
  $\noo R \rrm\,\le\,1$ and there is a subsequence $(l_k)_1^n \subset
   \{1,..,m\}$ such that $\noo \ds g^k
   \rrm_{\infty}\,=\,\bet<e_{l_k},\ds g^k>
  \rag$. From the rearrangement invariance of Y we deduce

\begin{eqnarray*}
 \noo \tau \rrm_Y &=& \noo \kla \bet\si_{l_k}<e_{l_k},g^k>\rag
     \mer_1^n \rrm_Y\\
 &\le& \noo \summ_1^n \si_{l_k}\,e_{l_k} \rrm_Y \\
 &\le& \noo \si\rrm_Y \: \le \:1 \pll.\\
\end{eqnarray*}

  Hence we obtain

\begin{eqnarray*}
  \noo \,\summ_1^n \noo TR\ds S(e_k) \rrm_F \, e_k \, \rrm_Z &=&
 \noo \,\summ_1^n  \kla \noo TJ(e_k) \rrm_F \,\noo\ds g^k \rrm_{\infty}
 \mer
      \, e_k \, \rrm_Z \\[+0.3cm]
   &\le& \pxe^n(T) \, \noo \tau \rrm_Y \: \le \: \pxe^n(T) \:
   .\\[-1.3cm]
\end{eqnarray*} \hfill $\Box$

\section{Gaussian cotype conditions}

\setcounter{lemma}{0}

As a consequence of Talagrand's factorization theorem for gaussian
processes
cotype conditions on $C(K)$-spaces can be reformulated with a quotient
formula. This was remarked by Pisier and Montgomery-Smith, see
\cite{MSM}.
We will give a prove for an arbitrary maximal, symmetric sequence
space.
Let us recall that $\lih$ is the space of sequences $\si\in\lif$ with

\[ \noo
\si\rrm_{\lih}\,:=\,\sup_{k\in\nz}\sqrt{\log(k+1)}\,\si_k^*\,<\,
 \infty . \]
\hz
\begin{lemma} \label{cotcon} Let X be a maximal, symmetric sequence
space,
   $T\in \com C(K),F)$ and $\nen$. Then we have for an absolut constant
   $c_1$

 \begin{eqnarray*}
   gc_X^n(T) &\sim_{c_1}& \sup\{\,\pi_{X,2}^n(TRD_{\si})\, |\, R \in
   \com \ce,
    E),\, \ds \in \com \ce,\ce)\,\mbox{with}\, \noo R\rrm\,,
    \noo \si \rrm_{\lih}\le 1 \, \} \, .\\
 \end{eqnarray*}
\end{lemma}
\hz
$\bf Proof:$ $\bf "\ge"$ W.l.o.g. we can assume that $\si_k\,=\,
  (\log(k+1))^{-1/2}$. Then it follows from \cite{LIP} that for all $u
  \in
  \com\lzn,\ce)$ we have

\[ \noo \summ_1^n g_k\,R\ds u(e_k)\rrm_{L_2(C(K))} \, \le \, c_1 \noo
R\rrm
   \, \noo u\rrm \pll . \]

With a glance on definition of $gc_X^n$ we see that the first
inequality
is proved. \hz

$\bf "\le"$ Let $(x_k)_1^n \in C(K)$ with

\[\noo \summ_1^n g_k x_k\rrm_{L_2(C(K))} \, \le \pll\, 1\pll .\]

By Talagrand's factorization theorem, see (*) in the preliminaries,
there are
$u \in \com \lzn,\ce)$ and $R \in \com \ce,C(K))$ with $\noo u\rrm \le
c_1$,
$\noo R \rrm \le1$ such that

 \[ RD_{\sigma}u(e_k) \: =\: x_k \pll, \]

and $\si_k\,=\,(\log(k+1))^{-1/2}$. Hence we deduce that

\begin{eqnarray*}
 \noo \summ_1^n \noo Tx_k\rrm_F \, e_k \rrm_X &=&
    \noo \summ_1^n \noo TR\ds u(e_k)\rrm_F \, e_k \rrm_X \\
  &\le& \pi_{X,2}^n(TR\ds) \, \noo u \rrm \\
 &\le& c_1\, \pi_{X,2}^n(TR\ds)\, \noo \summ_1^n g_k x_k
 \rrm_{L_2(C(K))}
   \pll .\\
\end{eqnarray*}
Taking the supremum over all sequences $(x_k)_1^n$ yields the
assertion.
\hfill $\Box$ \hz

Now we are able to prove the comparision theorem for gaussian and
Rademacher
cotype. \hz

\begin{theorem} \label{com} Let $2\!<\!q\!<\!\infty$ and X a q-convex
 maximal, symmetric sequence space. We set $Y\,=\,\dop(\lih,X)$. Then
 we
 have for all $T\in \com C(K),F)$ and $\nen$
\begin{enumerate}
\item   $ \pi_{Y,1}^n(T) \, \le \, rc_Y^n(T) \, \le \, \sqrt{2} \,
	 \pi_{Y,2}^n(T)\,  \le \, c_0 \, M^q(X) \kla
	 \frac{1}{2}\,-\,\frac{1}{q} \mer^{-1/q'}
    \, \pi_{Y,1}^n(T)\,. $
\item $gc_X^n(T) \pll \sim_{c_q} \pll rc_Y^n(T) \pll.$
\end{enumerate}
\end{theorem}
\hz

$\bf Proof:$ First we note that the q-convexity of X implies the
q-convexity
 of the maximal, symmetric sequence space Y. This can be seen exactly
 as in
 the proof of proposition\ref{qconvex} . Therefore the first assertion
 follows from theorem \ref{gema}, more precisely remark \ref{intthm},
 applied
 for Y. With the help of the previous Lemma \ref{cotcon}, applying
 theorem
 \ref{gema} for X and with the second quotient formula \ref{quot2} we
 obtain

 \begin{eqnarray*}
   gc_X^n(T) &\sim_{c_1}& \sup\{\,\pi_{X,2}^n(TRD_{\si})\, |\, R \in
   \com \ce,
    E),\, \ds \in \com \ce,\ce)\,\mbox{with}\, \noo R\rrm\,,
    \noo \si \rrm_{\lih}\le 1 \, \}\\
   &\sim_{c_q(X)}&  \sup\{\,\pi_{X,1}^n(TRD_{\si})\, |\, R \in \com
   \ce,
    E),\, \ds \in \com \ce,\ce)\,\mbox{with}\, \noo R\rrm\,,
    \noo \si \rrm_{\lih}\le 1 \, \}\\
   &=& \pi_{Y,1}^n(T) \pla . \\
 \end{eqnarray*}

Using the first assertion we see that the proof of the second assertion
is completed.
\hfill $\Box$
\hz

\begin{rem} \label{anwendung} Probably the most important applications
of
 the above theorem are given for gaussian cotype q and Rademacher
 cotype q
 operators when $q\!>\!2$.
\begin{enumerate}
\item In the case when $X\,=\,\ell_q$ it turns out that Y is in fact
 the Lorentz-Marcinkiewicz space $\ell_{q,q,-1/2}$. This space consists
 of all sequences $\si \in \lif$ such that

\[ \kla \summ_k \kla \frac{\si_k^*}{\sqrt{\log(k+1)}}  \mer^q
\,\mer^{1/q}
  \, < \, \infty \, . \]

\item If we want to calculate the cotype conditions for $(q,1)$-summing
 operators or Rademacher cotype q operators we have to solve the
 equation

 \[\ell_q \, = \, \dop(\lih,Y) \pll. \]

Again this is easy with the help of Lorentz-Marcinkiewicz spaces. The
space
$\ell_{q,q,-1/2}$ with the norm

\[ \noo \si \rrm_{\ell_{q,q,-1/2}} \,:=\, \kla \summ_k \kla
 \si_k^* \,\sqrt{\log(k+1)}\,  \mer^q \,\mer^{1/q} \]

solves the problem up to some constant. In order to apply theorem
\ref{com}
we have to check the r-convexity of $\ell_{q,q,-1/2}$ for some
$r\!>\!2$.
If we identify $\ell_{q,q,-1/2}$ with a space $\ell_{f,q}$ this easily
follows from remark \ref{example}. Indeed, f is given by

\[ f(t) \, := \, t^{1/q}\,\sqrt{\,\log(t+1)} \pl, \]

which satisfies $\be_f\, =\,\al_f\,=\, \frac{1}{q}$.
\end{enumerate}
\end{rem}\hz

In the following we will state further applications of theorem
\ref{com}.
\hz

\begin{cor} Let $2\!<\!q\!<\!\infty$ and X a q-convex maximal,
symmetric
sequence space then there is a constant $c$ depending on q and X only
such
that for all $\nen$ and $T\in \com C(K),F)$ with $rank(T)\,\le\,n$ one
has

 \[ gc_X(T) \,\le \, c \, gc_X^n(T) \pll.\]

\underline{Moreover}, the gaussian cotype constant is, up to $c$,
attained on
n disjoints functions in $C(K)$.
\end{cor} \hz

$\bf Proof:$ We set $Y = \dop(\lih,X)$. By theorem \ref{com} we have

\[ gc_X(T) \,\sim_c \, \pi_{Y,1}(T) \pll .\]

Therefore it remains to show that the $(Y,1)$-summing norm is attained
on
n vectors. Using Maurey's lemma about the extreme points of operators
from
$\lif^n$ to $C(K)$ (already used in the proof of proposition
\ref{quot2}),
see \cite{MAU}, it is then clear from that a restriction to n disjoint
blocs
is possible.

In theorem \ref{com} it was also observed that Y is q-convex. By
proposition \ref{qconvex} there is a maximal, symmetric sequence space
Z
with $Y\\cong\,\dop(Z,\ell_q)$. Furthermore, it is known that for the
computation of the $(q,2)$-summing norm of an operator with rank n only
n
vectors are needed, see for example \cite{DJ}.  Hence we can deduce
from
proposition \ref{quot1} and theorem \ref{gema}
\begin{eqnarray*}
 \pi_{Y,1}(T) &\le& \sup\{\,\pi_{q,2}(D_{\si}RT)\, |\, R \in \com F,
   \lif),\, \ds \in \com\lif,\lif),\,\mbox{with} \noo R\rrm,
   \noo \si \rrm_Z\le 1 \, \} \\
 &\le& \sqrt{2} \,\sup\{\,\pi_{q,2}^n(D_{\si}RT)\, |\, R \in \com F,
   \lif),\, \ds \in \com\lif,\lif),\,\mbox{with} \noo R\rrm,
   \noo \si \rrm_Z\le 1 \, \} \\
 &=& \sqrt{2} \, \pi_{Y,2}^n(T) \\
 &\le& \sqrt{2} \,c_q\, \pi_{Y,1}^n(T) \pl.\\[-1.3cm]
\end{eqnarray*}\hfill $\Box$\hz

In particular, the corollary works for $X=\ell_q$. For the so-called
"weak"
theory it is natural to replace $\ell_q$
by weak-$\ell_q$. More precisely, an operator $T \in \com E,F)$ is said
to be a $weak\, cotype\, q$ operator, if there exists a constant
$c\!>\!0$
such that for all $u\in \com \lzn,E)$ one has

\[ \sup_{k=1,..,n} \,k^{1/q}\, a_k(Tu) \pll \le \pll c \,\ell(u)
\pll.\]

The best constant c will be denoted by $\om c_q(T)$. It was essentially
remarked by Mascioni, see \cite{MAS}, that for $q\!>\!2$ another
definition
would have been possible. An operator $T\in \com E,F)$ is of weak
cotype q
if and only if there exists a constant $c>0$ such that

\[ \sup_{k \in \nz} k^{1/q} \noo Tx_k \rrm_F \pll \le \pll
     c \, \noo \summ_k g_k x_k \rrm_{L_2(E)} \]

for each sequence $(x_k)_k \subset E$ such that $\noo Tx_k \rrm$ is
non-increasing (for further information see also \cite{DJ1}). The next
proposition gives a characterization of weak cotype operators on
$C(K)$-spaces in terms of Weyl numbers.
\hz

\begin{cor} Let $2\!<q\!<\!\infty$. An operator $T \in \com C(K),F)$ is
 of weak cotype q if and only if

\[ \sup_{k\in\nz} \frac{k^{1/q}}{\sqrt{\log(k+1)}} \pl x_k(T) \pll
  < \pll \infty \pll .\]
\end{cor} \hz

$\bf Proof:$ By remark \ref{example} the space $X\,:= \,
\ell_{q,\infty} \,
 := \, \ell_{f,\infty}$ with $f(t)=t^{1/q}$ is r-convex for all
 $2\!<r\!<\!q$.
We observe that $Y\,:=\,\dop(\lih,X)$ coincides with $\ell_{g,\infty}$
where
$g(t)= t^{1/q}/\sqrt{\log(t+1)}$. Using Mascioni's observation above we
deduce
from theorem \ref{com} that T is of weak cotype q if and only if T is
$(Y,2)$-summing.

If T is $(Y,2)$-summing and $u \in \com \ell_2,C(K))$ we can apply a
lemma due
to Lewis, see \cite{PIE}, which guarantees for all $\vare\!>\!0$ the
existence
of an orthonormal system $(o_k)_k \subset \ell_2$ with $\kla \noo
Tu(o_k)
\rrm_F\mer_k $ decreasing and

\[ a_k(Tu) \pll \le \pll (1+\vare) \pl\noo Tu(o_k) \rrm_F \pll .\]

Therefore we deduce

\begin{eqnarray*}
 \sup_{k\in\nz} \frac{k^{1/q}}{\sqrt{\log(k+1)}} \pl a_k(Tu)
 &\le& (1+\vare) \pl
 \sup_{k\in\nz} \frac{k^{1/q}}{\sqrt{\log(k+1)}} \pl \noo Tu(o_k)
 \rrm_F \\
 &\le& (1+\vare) \pl\pi_{Y,2}(T) \pl \om_2(u(o_k))_k \\
 &\le& (1+\vare) \pl \pi_{Y,2}(T)\pl \noo u\rrm \pll .\\
\end{eqnarray*}

Taking the infimum over all $\vare$ and the supremum over all $u \in
\com \ell_2,C(K))$ with norm less than 1 we obtain

\[ \sup_{k\in\nz} \frac{k^{1/q}}{\sqrt{\log(k+1)}} \pl x_k(T) \pll
  \le  \pll \pi_{Y,2}(T) \pll .\]

Vice versa, let us assume that the sequence of Weyl numbers is in Y.
Let
 $(x_k)_k \in C(K)$ with $\om_2(x_k)_k\le 1$. There is no restriction
to assume that $\noo Tx_k\rrm_F$ is decreasing. If we define $u_n
\,:=\,
\summ_1^n e_k \otimes x_k$ we can deduce from an inequality of K\"onig,
see
\cite{PIE},

\begin{eqnarray*}
 n^{1/2} \, \noo Tx_n\rrm &\le& \kla \summ_1^n \noo Tx_k \rrm^2
 \mer^{1/2}
   \pll \le \pll \pi_2(Tu_n) \\
 &\le& c_1 \, \summ_1^n \frac{a_k(Tu_n)}{\sqrt{k}} \\
 &\le& c_1 \, \summ_1^n \frac{ (\log(k+1))^{1/2} }{k^{1/2+1/q}} \pll
   \noo\summ_1^n x_k(T)\, e_k\rrm_Y \pll \noo u \rrm \\
 &\le& c_1 \sqrt{\log(n+1)} \pll \frac{n^{1/2-1/q}}{1/2-1/q}
   \pll \noo \summ_k x_k(T)\, e_k\rrm_Y \pl .\\
\end{eqnarray*}

Taking the supremum over all $\nen$ we have shown that T is
$(Y,2)$-summing.
\hfill $\Box$


\hz
1991 Mathematics Subject Classification: 47B38, 47A10, 46B07.
\begin{quote}
Marius Junge

Mathematisches Seminar der Universit\H{a}t Kiel

Ludewig-Meyn-Str. 4

W-2300 Kiel 1

Germany
\end{quote}


\begin{thebibliography}{10001}
\bibitem[COB]{COB} F. Cobos: {\sl On the Lorentz-Marcinkiewicz Operator
ideal.}: Math. Nachr. 126.(1986), 281-300.
\bibitem[DJ]{DJ} M. Defant and M.Junge: {\sl On absolutely summing
operators with apllication to the (p,q)-summing norm with few vectors};
J. of Functional Ana. 103 (1992), 62-73.
\bibitem[DJ1]{DJ1} M. Defant and M. Junge: {\sl Random variables in
weak type p spaces}; Arch. Math. 58 (1992), 399-406.
\bibitem[MAS]{MAS} V. Mascioni: {\sl On weak cotype and weak type in
Banach spaces}; Note di Matematica Vol VIII-n.1(1988), 67-110.
\bibitem[MAU]{MAU} B. Maurey: {\sl Type et cotype dans les espaces
munis d'un structure localement inconditionelle}; S\'{e}minaire
Maurey-Schwartz 73-74, Ecol- Polyt., Exp. no. 24-25.
\bibitem[MSM]{MSM} S.J.Montgomery-Smith: {\sl The Gaussian cotype of
operators from $C(K)$}; Isr. J. Math. 68 (1989), 123 - 128.
\bibitem[LET]{LET} M. Ledoux and L.Talagrand: {\sl Probability in
Banach spaces}. Berlin Heidelberg New York: Springer 1991.
\bibitem[LIP]{LIP} W. Linde and A. Pietsch: {\sl Mappings of gaussian
cylindrical measures in Banach spaces}; Theory Probab. Appl. 19 (1974),
445-460.
\bibitem[LTI]{LTI} J. Lindenstrauss and L. Tzafriri: {\sl Classical
Banach spaces I, sequence spaces}; Springer Berlin Heidelberg New York
1977.
\bibitem[LTII]{LTII} J. Lindenstrauss and L. Tzafriri: {\sl Classical
Banach spaces II, function  spaces}; Springer Berlin Heidelberg New
York 1979.
\bibitem[PIE]{PIE} A. Pietsch: {\sl Eigenvalues and s-numbers of
operators}; Cambridge University Press, 1987.
\bibitem[TAL]{TAL} M. Talagrand: {\sl Cotype of operators from $C(K)$};
Invent. math. 107 (1992), 1-40.
\bibitem[TA1]{TA1} M. Talagrand: {\sl Regularity of Gaussian
processes}; Acta. Math. 159 (1987), 99-149.
\bibitem[TJM]{TJM} N. Tomczak-Jaegermann: {\sl Banach-Mazur distances
and finite-dimensional operator ideals}; Longmann, 1988.
\end{thebibliography}
\end{document}